\documentclass[secthm]{elsart}

\usepackage{amssymb}
\usepackage{amsmath}
\usepackage{longtable}

\def\af#1{\mathbb A^{#1}}
\def\al{\alpha}
\def\as#1{\renewcommand\arraystretch{#1}}
\def\bbb{\mathbf{B}}

\def\cc{{\mathcal C}}
\def\comb#1#2{\as{1}\left(\!\begin{array}{c}
#1\\#2
\end{array}\!\right)}
\def\combr#1#2{\as{1}\left(\!\!\left(\!\begin{array}{c}
#1\\#2
\end{array}\!\right)\!\!\right)}

\def\dg{\op{deg}}
\def\di{\op{diag}}
\def\dm{\op{dim}}
\def\e{\medskip}
\def\eps{\epsilon}
\def\exp{\op{exp}}
\def\ff#1{\mathbb F_{q^{#1}}}
\def\fix#1#2{\operatorname{Fix}_{\gamma}(#1,#2)}
\def\fixr#1#2{\operatorname{Fix}_{\gamma}(\!(#1,#2)\!)}

\def\fq{\mathbb{F}_q}
\def\g{\Gamma}
\def\ga{\gamma}
\def\gal{\op{Gal}}
\def\gen#1{\big\langle\, {#1} \,\big\rangle}
\def\gl#1#2{\op{GL}_{#1}(#2)}
\def\ig{\ll_{\ga}^{\op{pr}}}

\def\igk{\ll_{\ga}^{\op{pr}}(k)}

\def\igg{\ll_{\ga,G}^{\op{pr}}}
\def\imp{\,\Longrightarrow\,}
\def\iso{\,\stackrel{\mbox{\tiny $\sim\,$}}{\lra}\,}
\def\k{\op{Ker}}
\def\kb{\overline{k}}
\def\la{\lambda}

\def\ll{{\mathcal L}}
\def\lra{\longrightarrow}
\def\lcm{\op{lcm}}
\def\md#1{\ \mbox{\rm(mod }{#1})}
\def\N{\mathbb{N}}

\def\ni{\noindent}

\def\op{\operatorname}
\def\ord{\op{ord}}

\def\pgl{\operatorname{PGL}_{N+1}(k)}
\def\pr#1{\mathbb P^{#1}}
\def\prn{\mathbb P^{N}(k)}
\def\Q{\mathbb{Q}}

\def\sg{\sigma}
\def\sii{\,\Longleftrightarrow\,}
\def\st{\op{st}}
\def\tq{\,\,|\,\,}

\def\xg{\op{Fix}_{\ga}(X)}
\def\Z{\mathbb{Z}}

\begin{document}

\begin{frontmatter}



\title{Orbits of rational $n$-sets of projective spaces under the action of the linear group}


\author{Ricard Mart\'i},
\author{Enric Nart\corauthref{corr}\thanksref{asse}}

\address{Universitat Aut\`onoma de Barcelona, Departament de Matem\`atiques \\ 08193 Bellaterra, Barcelona, Spain}
\ead{nart@mat.uab.cat\quad Fax number: +34935812790}

\corauth[corr]{Corresponding author}
\thanks[asse]{Supported by the project MTM2006-11391 from the Spanish MEC}

\begin{abstract}
For a fixed dimension $N$ we compute the generating function of the numbers $t_N(n)$ (respectively $\overline{t}_N(n)$) of $\pgl$-orbits of rational $n$-sets
 (respectively rational $n$-multisets) of the projective space $\pr{N}$ over a finite field $k=\fq$. For $N=1,2$ these results provide concrete formulas for $t_N(n)$ and $\overline{t}_N(n)$ as a polynomial in $q$ with integer coefficients.
\end{abstract}


\end{frontmatter}

\section*{Introduction}
There are many examples of deep properties of geometric objects relying on combinatorial properties of unordered structures. In particular, $n$-sets of projective spaces have deserved the attention of geometers since a long time. The work of Coble at the beginning of the last century is an outstanding contribution to the study of geometric structures underlying $n$-sets of projective spaces.
A revision of this work in modern language can be found in the book \cite{dol} of Dolgachev and Ortland. 

If we are interested in arithmetic properties of the geometric objects associated to $n$-sets, we need to consider {\it rational} $n$-sets; that is, $n$-sets defined over the ground fied $k$ we are interested in. Let $\kb$ be a fixed algebraic closure of the field $k$; an $n$-set $S=\{P_1,\dots,P_n\}$ of a projective space $\pr{N}(\kb)$ is $k$-rational if $S$ is invariant under the action of the absolute Galois group
$$
\{\sg(P_1),\dots,\sg(P_n)\}=\{P_1,\dots,P_n\},\quad\forall\sg\in\gal(\kb/k).
$$  
Thus, $S$ is the disjoint union of orbits of points of $\pr{N}(\kb)$ under the action of the Galois group. We denote by 
$$
\comb{\pr{N}}{n}(k):=\comb{\pr{N}(\kb)}{n}^{\gal(\kb/k)},\quad 
\combr{\pr{N}}{n}(k):=\combr{\pr{N}(\kb)}{n}^{\gal(\kb/k)},
$$
the respective sets of $k$-rational $n$-sets and $n$-multisets of $\pr{N}(\kb)$. 

If we apply to a rational $n$-set a $k$-automorphism of $\pr{N}$ we obtain an equivalent $n$-set, in the sense that the underlying geometric objects of both $n$-sets will have the same geometric and arithmetic properties. 

The aim of this paper is the computation of the number of $\pgl$-orbits of rational $n$-sets and $n$-multisets of projective spaces $\pr{N}$ defined over a finite field $k$. That is, we want to find closed formulas for the numbers:
$$
t_N(n):=\left|\pgl\backslash \comb{\pr{N}}{n}(k)\right|,\quad
\overline{t}_N(n):=\left|\pgl\backslash \combr{\pr{N}}{n}(k)\right|.
$$ 

There is an extensive literature on the enumeration of orbits of {\it pointwise rational} $n$-sets; that is, $n$-sets $S=\{P_1,\dots,P_n\}$ such that all $P_i$
are $k$-rational points (have homogeneous coordinates in $k$). This is due to the fact that these orbits are in correspondence with isometry classes of linear codes \cite{fr},
\cite[Sec. 3.2]{bfkwz}, \cite{ln}, \cite{mn}, \cite{mn2}. However, to our knowledge the enumeration of rational $n$-sets has not been considered so far, with the exception of \cite{lmnx}, where the numbers $t_1(n)$ were computed.

Let us illustrate both the role of global (not pointwise) rationality and the action of the linear group with an example. It is well-known that the hyperelliptic curves over an algebraically closed field (of zero or odd characteristic) are parameterized by $n$-sets of $\pr1$; if $S=\{P_1,\dots,P_n\}$ is a $n$-set of $\pr1$ and we attach to each $P_i$ an affine coordinate $a_i$, we can consider the hyperelliptic curve given by the Weierstrass equation $$y^2=(x-a_1)(x-a_2)\cdots(x-a_n).$$ If we want to classify hyperelliptic curves defined over a non-algebraically closed field $k$ we are led to consider $k$-rational $n$-sets of $\pr1$. If $S$ is a pointwise rational $n$-set the above construction provides a curve with all Weierstrass points defined over $k$. These curves are a small part of the family of hyperelliptic curves defined over $k$, given by Weierstrass equations $y^2=f(x)$, with $f(x)$ an arbitrary separable  polynomial with coefficients in $k$. Finally, it is easy to check that rational $n$-sets in the same orbit by the action of $\op{PGL}_2(k)$ determine $k$-isomorphic curves. Thus, the numbers $t_1(n)$ count essentially $k$-isomorphism classes of hyperelliptic curves defined over $k$ \cite{lmnx}.

Another interesting example is given by the $7$-sets of the projective plane. Over an algebraically closed field certain $\op{PGL}_3(k)$-orbits of $7$-sets classify non-hyperelliptic curves of genus three with a fixed $2$-level structure \cite[Chap. IX]{dol}. Thus, our enumeration results may provide information on the number of $k$-rational points of the moduli space of such objects, and from these numbers one can derive further geometric and arithmetic information on this space.

Our main result is a computation of the generating function of the $t_N(n)$, $\overline{t}_N(n)$ for fixed $N$ (Theorem  \ref{ffinal}). It is well-known that the generating function of the number of orbits of pointwise rational $n$-sets of $\prn$ can be expressed in terms of the cycle index of P\'olya \cite[3.2.16]{bfkwz}. This cycle index is a polynomial in several variables that carries all information about the lengths of the cycles of all elements of $\pgl$ acting as permutations of $\prn$. Thus, this instrument is not able to provide information on rational $n$-sets. However, in \cite{mn2} we have found a refinement of the cycle index that leads to more effective formulas for the generating function of the number of orbits of pointwise rational $n$-sets of $\prn$. In section \ref{secquatre} we extend these ideas to rational $n$-sets.        
We introduce certain equivalence relation in the set $\cc$ of conjugacy classes of $\g:=\pgl$; if $\mathcal{S}$ is the quotient set, for each class $\al\in \mathcal{S}$  (we call $\al$ a {\it subtype}) we construct a poset $\ll(\al)$ of certain $\ga$-invariant linear subvarieties of $\pr{N}(\kb)$ classified under the action of the Galois group. The nodes $V$ of this poset $\ll(\al)$ carry a weight $(\dm V,\exp V,\dg V)$ of three numerical values {\it dimension}, {\it exponent} and {\it degree}. We count the elements in the same subtype $\al$ in the weighted sum: $$M_{\al}:=\sum_{\ga\in \al}\vert\g_\ga\vert^{-1},$$ where $\g_{\ga}$ is the centralizer of $\ga$ in $\g$. Then, we consider a {\it $G$-exponent index} ($G$ stands for ``Galois"): 
$$
\op{E}_G(\op{PGL}_{N+1},\pr{N}):=\sum_{\al\in\mathcal{S}}M_{\al}\prod_{V\in\ll(\al)}z_{\al,V}\in\Q[\{z_{\al,V}\}],
$$
able to express the generating function of the numbers $t_N(n)$, $\overline{t}_N(n)$ as:
\begin{equation}\label{main}\as{1.6}
\begin{array}{l}
\sum_{n\in \N}t_N(n)x^n=\op{E}_G(\op{PGL}_{N+1},\pr{N})\left[z_{\al,V}=f_{\al,V}(x)\right], \\ \sum_{n\in
\N}\overline{t}_N(n)x^n=\op{E}_G(\op{PGL}_{N+1},\pr{N})\left[z_{\al,V}=\overline{f}_{\al,V}(x)\right],
\end{array}
\end{equation}
for certain explicitly given functions $f_{\al,V}(x)$, $\overline{f}_{\al,V}(x)$ (cf.Theorem \ref{ffinal}). 

Section \ref{secu} is devoted to the computation of the total number of rational $n$-sets of an arbitrary quasiprojective variety $V$; the generating function of these numbers is easily expressed in terms of the zeta function of $V$. In section \ref{secdos} we prove a crucial step in the enumeration of $\pgl$-orbits of rational $n$-sets of $\pr{N}$: for any $\ga\in\pgl$, the quotient variety $\pr{N}/\gen{\ga}$ has the same zeta function than $\pr{N}$ (Theorem \ref{tdos}). In section \ref{sectres} we introduce the concept of {\it proper} $\ga$-invariant linear subvariety of $\pr{N}(\kb)$ for each $\ga\in\pgl$; this is the main ingredient in the construction of the set $\mathcal{S}$ of subtypes and their associated posets $\ll(\al)$.  
Finally, in section \ref{seccinc} we restrict our attention to the cases $N=1,\,2$ and we carry out an explicit computation of all the ingredients of (\ref{main}) in terms of combinatorial data independent of the group structure of $\pgl$ and the action of its elements as permutations of $\pr{N}(\kb)$. This allows one to obtain explicit expressions for the numbers $t_N(n)$, $\overline{t}_N(n)$ as polynomials with integer coefficients in  the cardinality $q$ of the ground field.

\ni{\bf Conventions and notation. }Throughout the paper we fix a finite field $k=\fq$ of characteristic $p$ and an algebraic closure $\kb$ of $k$.
For any integer $r\ge1$ we denote by $k_r=\ff{r}$ the unique
extension of degree $r$ of $k$ in $\kb$.
We denote by $\sigma(x)=x^q$ the $q$-Frobenius automorphism of $\kb$, which is a topological generator of the absolute Galois group  $\gal(\kb/k)$ as a profinite group.

\section{Rational $n$-sets of quasiprojective varieties}\label{secu}
Let $V$ be a quasiprojective algebraic variety  defined over $k$. The {\it variety of $n$-multisets of $V$} is by definition the symmetric product of $V$ with itself $n$ times. The {\it variety of $n$-sets of $V$} is the nonsingular locus of the former variety; that is, the open subvariety formed by the unordered $n$-tuples of points of $V$ without repetitions. We denote these varieties by:    
$$
\comb{V}{n}\subseteq \combr{V}{n}:=(V\times\stackrel{n)}\cdots\times V)/S_n.
$$
The $k$-rational $n$-sets and $n$-multisets of $V$ are respective $k$-rational points of these varieties; that is, elements of the sets 
$$\comb{V}{n}(k) := \comb{V(\kb)}{n}^{\gal(\kb/k)},\quad\combr{V}{n}(k) := \combr{V(\kb)}{n}^{\gal(\kb/k)}.
$$

In this section we want to compute the total number of rational $n$-sets and $n$-multisets of $V$
$$a_{V}(n):=\left|\comb{V}{n}(k)\right|,\qquad \bar{a}_{V}(n):=\left|\combr{V}{n}(k)\right|.
$$
By convention, $a_{V}(0)=1=\bar{a}_{V}(0)$. Our first step is to express $a_{V}(n)$, $\bar{a}_{V}(n)$ in terms of the numbers of orbits of given length of $V(\kb)$ under the action of $\gal(\kb/k)$.
\begin{defn}
For any $P \in V(\kb)$, we denote by $O_{\sg}(P)$ the orbit of $P$ under the action of $\gal(\kb/k)$. We call $O_{\sg}(P)$ the $\sg$-orbit of $P$.

The degree of $P$ is $\dg(P):=|O_{\sg}(P)|$; in other words, $\dg(P)$ is the minimum positive integer $r$ such that $P$ lies in
$V(k_r)$.

Finally, we denote by $b_r$ the number of $\sg$-orbits of length $r$ of points of $V(\kb)$:
$$
b_r:=\left|\{O_{\sg}(P)\tq \dg(P)=r\}\right|=\frac 1r\left| \{P\in V(\kb)\tq
\dg(P)=r\}\right|.
$$
\end{defn}
If we think a $k$-rational  $n$-set or $n$-multiset of $V$ as a disjoint union of $\sg$-orbits of different length
we can express $a_{V}(n),\,\bar{a}_{V}(n)$ in terms of $b_1,b_2,\dots,b_n$:
$$\as{2.2}
\begin{array}{l}
a_{V}(n) = \sum_{s_1+2s_2+\cdots+ns_n=n}\comb{b_1}{s_1}\comb{b_2}{s_2}\dots\comb{b_n}{s_n},
\\
\bar{a}_{V}(n) = \sum_{s_1+2s_2+\cdots+ns_n=n}\combr{b_1}{s_1}\combr{b_2}{s_2}\dots\combr{b_n}{s_n},
\end{array}
$$
where $s_i$ is the number of $\sg$-orbits of degree $i$ in each
$n$-set or $n$-multiset, and we understand that $\comb{b_i}{s_i}=0$ if $s_i>b_i$. These expressions are bad explicit formulas
for $a_{V}(n),\,\bar{a}_{V}(n)$ because their evaluation is extremely inefficient; however, they lead to a computation of the generating function of these numbers:
\begin{equation}\label{fbes}
\as{2.}
\begin{array}{l}
f_{V}(x):=\sum_{n\ge0}a_{V}(n)x^n=\prod_{r\ge 1}(1+x^r)^{b_r},\\
\bar{f}_{V}(x):=\sum_{n\ge0}\bar{a}_{V}(n)x^n=\prod_{r\ge 1}(1-x^r)^{-b_r},
\end{array}\end{equation}
closely related to the zeta function of $V$ over $k$:
$$\op{Z}(V/k,x)=\exp\left(\sum_{r\ge1}\frac1{r}N_rx^r\right),\qquad N_r=N_r(V):=|V(k_r)|.$$
 The families $\{N_r\}_{r\ge 1}$ and $\{b_r\}_{r\ge1}$ are determined one by each other and their relationship is synthesized
in the following expression of the zeta function as an infinite product:
$$\op{Z}(V/k,x)=\prod_{r\ge1}(1-x^r)^{-b_r}.$$

From this expression and (\ref{fbes}) we deduce the main result of this section:
\begin{thm}\label{tu}
For any quasiprojective variety $V$ defined over $k$:
$$
f_{V}(x)=\op{Z}(V/k,x)/\op{Z}(V/k,x^2),\qquad \bar{f}_{V}(x)=\op{Z}(V/k,x).
$$
\end{thm}
\begin{cor}\label{multipl}
Let $W\subseteq V$ be any subvariety of $V$, which is also defined over $k$, and let
 $U=V\setminus W$ be the complementary subvariety. Then $f_{V}(x)=f_{W}(x)f_{U}(x)$.
\end{cor}
\begin{pf}
Clearly $\ \op{Z}(V/k,x)=\op{Z}(W/k,x)\op{Z}(U/k,x)$. \hfill{$\Box$}
\end{pf}
In the rest of the section we apply Theorem \ref{tu} to obtain explicit formulas for $a_{V}(n)$ for different particular varieties $V$. These computations will be used in section \ref{seccinc} to obtain explicit formulas for $t_2(n)$.
 
\subsection{Explicit formulas for $a_{V}(n)$: subvarieties of $\af{N}$}
Since $\op{Z}(\af{N}/k,x)=(1-q^Nx)^{-1}$, we get immediately from Theorem \ref{tu}:
\begin{equation}\label{fafn}
f_{\af{N}}(x)=\frac{1-q^Nx^2}{1-q^Nx}=q^{-N}+x+\frac{1-q^{-N}}{1-q^Nx}=q^{-N}+x+(1-q^{-N})\sum_{n\ge0}q^{Nn}x^n. 
\end{equation}

Thus, we get a closed formula for the $n$-coefficient of this series: 

\begin{prop}\label{3}
For all $N \ge0$,
$\as{1.4}\quad
a_{\af{N}}(n)=\left\{\begin{array}{ll}
q^{nN},&\mbox{ if }n\le 1,\\
q^{nN}-q^{(n-1)N},&\mbox{ if }n\ge2.
\end{array}    \right.\as{1}$
\end{prop}

For $N=1$ the number $a_{\af{1}}(n)$ counts the number of monic separable polynomials 
of degree $n$ with coefficients in $\fq$. This formula has been rediscovered several times in the literature (cf. \cite{fjk}, \cite{bg}, \cite{lmnx}). To our knowledge the formula for $a_{\af{N}}(n)$, $N > 1$, is new.

We need similar formulas for certain open subvarieties of $\af{N}$. If $\{\ast\}\subseteq \af{N}$ is the $0$-dimensional variety given by a $k$-rational point and
$L,L'\subseteq \af2$ are two non parallel lines, we get by Corollary \ref{multipl}
$$
f_{\af{N}\setminus\{\ast\}}(x)=\frac{f_{\af1}(x)}{1+x},\quad 
f_{\af2\setminus
L}(x)=\frac{f_{\af2}(x)}{f_{\af1}(x)}, \quad
f_{\af{2}\setminus (L\cup L')}(x)=\frac{f_{\af2}(x)(1+x)}{f_{\af1}(x)^2},
$$
since $f_{\{\ast\}}(x)=1+x$, $f_L(x)=f_{\af1}(x)$ and $f_{L\cup L'}=f_{\af1}(x)^2/f_{\{\ast\}}(x)$. We obtain explicit formulas for   
the number of rational $n$-sets in these varieties just by developping these rational functions as a series and finding the $n$-th coefficient. 
 
\begin{prop}\label{lines}
$\quad a_{\af{N}\setminus\{\ast\}}(n)=(q^N-1)\dfrac{q^{nN}-(-1)^n}{q^N+1},\quad
\forall n\ge 1$.
$$
\as{2.} a_{\af2\setminus L}(n)=\left\{
\begin{array}{ll}
\dfrac{q^2-1}{q^2+q+1}\left(q^{2n}-q^{n/2}\right),&\mbox{
if $\,0<n$ even}, \\
\dfrac{q-1}{q^2+q+1}\left((q+1)q^{2n}+q^{(n+1)/2}\right),&\mbox{
if $\,n$ odd}.
\end{array}
\right.
$$
Moreover, for $n$ even, $n>0$:
$$
a_{\af{2}\setminus (L\cup L')}(n)=\dfrac{q^4-1}{(q^2+q+1)^2}\left(q^{2n}-q^{n/2}\left(\dfrac
n2\,\dfrac{(q^3-1)(q-1)}{q^4-1}+1\right)\right), 
$$
whereas for $n$ odd, $n>1$, the value of $a_{\af{2}\setminus (L\cup L')}(n)$ is:
$$ 
\dfrac{q^4-1}{(q^2+q+1)^2}\left(q^{2n}+q^{(n-1)/2}
\left(\dfrac{n-1}2\,\dfrac{q^3-1}{q^2+1}-\dfrac{(q-1)(2q^2+q+1)}{q^4-1}\right)\!\right).
$$
\end{prop}

\subsection{Explicit formulas for $a_{V}(n)$: subvarieties of $\pr{N}$}
By the usual stratification of $\pr{N}$ as a union of affine spaces we get from (\ref{fafn}):
\begin{equation}\label{fpr}
f_{\pr{N}}(x)=\frac{(1-x^2)(1-qx^2)\cdots(1-q^Nx^2)}{(1-x)(1-qx)\cdots(1-q^Nx)}.
\end{equation}
For small values of $N$ we can find easily a closed expression for the $n$-th coefficient of this series. 
\begin{prop}\label{pset}
$$\as{1.2}
a_{\pr{1}}(n)=\left\{
\begin{array}{ll}
q+1,&\mbox{ if } \,n=1,\\
q^2,&\mbox{ if } \,n=2,\\
q^n-q^{n-2},&\mbox{ if } \,n\ge3
\end{array}\right.
$$$$\as{1.2}
a_{\pr{2}}(n)=\left\{
\begin{array}{ll}
q^2+q+1,&\mbox{ if }\,n=1,\\
q^4+q^3+q^2,&\mbox{ if }\,n=2,\\
q^6+q^5+q^4-q^2-q,&\mbox{ if }\,n=3,\\
(q^6+q^5+q^4-q^2-q-1)q^{2n-6},&\mbox{ if }\,n\ge4.
\end{array}\right.
$$
\end{prop}
Althoug we are not going to use it in section \ref{seccinc}, let us display the computation of $a_{\pr{N}}(n)$ for $n>N+1$ and arbitrary $N$. We can write (\ref{fpr}) in the form
$$
f_{\pr{N}}(x)=\frac{\prod_{0\le i\le\lfloor\frac N2\rfloor}(1+q^ix)\prod_{1\leq i\leq N,\,
i\, \mbox{\tiny odd}}(1-q^ix^2)}{\prod_{\frac N2<i\le
N}(1-q^ix)}=a(x)+\frac{b(x)}{c(x)},
$$with polynomials $a(x),\,b(x),\,c(x)$ such that $\deg a(x)=N+1$, $\deg b(x)<\deg c(x)$.  The coefficients of $a(x)$ distort the values of $a_{\pr{N}}(n)$, for $n\leq N+1$, but for $n> N+1$ they depend only on 
the fraction $b(x)/c(x)$, which after decomposition into a sum of elementary fractions leads to 
\begin{equation}\label{eonze}
a_{\pr{N}}(n)=\sum_{\frac N2<i\leq N}\lambda_iq^{in},\quad \forall
n>N+1.
\end{equation}

\begin{prop}
For any positive integer $r$ let
$$
\Lambda(r):=(1-q^{-1})(1-q^{-2})\cdots(1-q^{-r})=q^{-r(r+1)/2}(q-1)(q^2-1)\cdots(q^r-1),
$$and take $\Lambda(0)=1$ by convention. Then,
$$
a_{\pr{N}}(n)=\sum_{\frac N2<i\leq
N}\dfrac{(-1)^{N-i}\Lambda(2i)}{\Lambda(2i-N-1)\Lambda(i)\Lambda(N-i)}\,
q^{in-(N-i)(N-i+1)/2},\quad \forall n>N+1.
$$
\end{prop}
\begin{pf}
Let $P(x)=(1-x)(1-qx)\cdots(1-q^Nx)$. We want to find constants $\lambda_i$ uniquely determined by:
$$
\frac{P(x^2)}{P(x)}\in\sum_{i=0}^N \frac{\lambda_i}{1-q^ix}+\Z[x].
$$
This relationship is equivalent to:
$$
P(x^2)\,\in\, \sum_{i=0}^N\, \lambda_i\,\frac{P(x)}{1-q^ix}\,+\,P(x)\,\Z[x],
$$
and we can isolate $\la_i$ by taking $x=q^{-i}$:
$$
\lambda_i=\frac{P(q^{-2i})}{\left[P(x)/(1-q^ix)\right]_{x=q^{-i}}}.
$$
We get $\lambda_i=0$ for $i=0, 1,\dots,\lfloor\frac N2\rfloor$, and 
$$
\lambda_i=\left(\prod_{r=0}^N (1-q^{r-2i})\right)\Big/\left(\prod_{s=0,\,s\ne i}^N
(1-q^{s-i})\right), \quad \frac N2<i\le N.
$$
The numerator of this fraction is $\Lambda(2i)\big/\Lambda(2i-N-1)$, and we can compute separatedly the factors with $s<i$ and $s>i$
of the denominator
$$
\prod_{s=0}^{i-1}(1-q^{s-i})=\Lambda(i),\qquad
\prod_{s=i+1}^N(1-q^{s-i})=(-1)^{N-i}\Lambda(N-i)q^{(N-i)(N-i+1)/2}.
$$ The proof ends by introducing these expressions in (\ref{eonze}).\hfill{$\Box$}
\end{pf}

Finally, we need to compute $a_{V}(n)$ in some special cases:
$$V_1=\pr1\setminus O_{\sg}(P),\quad V_2=\pr2\setminus \{\ast\},\quad V_3=\pr2\setminus O_{\sg}(Q),$$$$V_4=\pr2\setminus \{P_1,P_2,P_3\},\quad V_5=\pr2\setminus \{P_1,O_{\sg}(P)\}$$ 
where $P$ denotes a point of degree two, $Q$ a point of degree three and $P_1,\,P_2,\,P_3$ points of 
degree one. By Corollary \ref{multipl}, the respective generating functions are:
$$f_{V_1}(x)=f_{\pr1}(x)/(1+x^2),\quad f_{V_2}(x)=f_{\pr2}(x)/(1+x),\quad f_{V_3}(x)=f_{\pr2}(x)/(1+x^3),$$$$ f_{V_4}(x)=f_{\pr2}(x)/(1+x)^3,\quad f_{V_5}(x)=f_{\pr2}(x)/(1+x)(1+x^2).$$ 

\begin{prop}\label{qq}For $n>0$:
$$\as{2.}
a_{V_1}(n)=\left\{
\begin{array}{ll}
\dfrac{q+1}{q^2+1}\,\left((q-1)q^n+(-1)^{(n-1)/2}(q+1)\right),&\mbox{ if $\,n$ odd},\\
\dfrac{q^2-1}{q^2+1}\,(q^n-(-1)^{n/2}),&\mbox{ if $\,n$ even}.
\end{array}\right.
$$
$$a_{V_2}=
\left\{
\begin{array}{ll}
q^2+q,&\mbox{ if }\,n=1,\\
q^4+q^3-q,&\mbox{ if }\,n=2,\\
(q^4+q^3-q-1)q^{2n-4},&\mbox{ if }\,n\ge3.
\end{array}\right.
$$
$$\as{2.}
a_{V_3}(n)=\left\{
\begin{array}{ll}
\dfrac{q^2+q+1}{q^4-q^2+1}\,((q^2-1)q^{2n}+(-1)^{(n-1)/3}),&\mbox{ if }n\equiv 1\md3,\\
\dfrac{q^2+q+1}{q^4-q^2+1}\,((q^2-1)q^{2n}+(-1)^{(n-2)/3}q^2),&\mbox{ if }n\equiv 2\md3,\\
\dfrac{(q^3-1)(q+1)}{q^4-q^2+1}\,(q^{2n}-(-1)^{n/3}),&\mbox{ if }n\equiv 0\md3.
\end{array}\right.
$$
\begin{multline*}
a_{V_4}(n)=\dfrac{q-1}{(q^2+1)^2}\left[(q^3+2q^2+2q+1)q^{2n}+\right.\\\left.+(-1)^n\left((n-1)q^3-(n+2)q^2+(n-2)q-(n+1)\right)\right].
\end{multline*}
$$\as{2}
a_{V_5}(n)=\left\{
\begin{array}{ll}
\dfrac{q+1}{q^4+1}\,((q^3-1)q^{2n}+(-1)^{(n-1)/2}q(q+1)),&\mbox{ if $\,n$ odd},\\
\dfrac{(q^3-1)(q+1)}{q^4+1}\,\left(q^{2n}-(-1)^{n/2})\right),&\mbox{ if $\,n$ even}.
\end{array}\right.
$$
\end{prop}

\section{Zeta function of the quotient of $\pr{N}$ by an automorphism}\label{secdos}
The aim of this section is to prove the following result: 

\begin{thm}\label{tdos}
For any $\ga\in \pgl$, let $\pr{N}/\ga$ be the quotient variety of $\pr{N}$ by the finite cyclic group generated by $\ga$. Then, 
$\op{Z}\left((\pr{N}/\ga)/k,x\right)=\op{Z}(\pr{N}/k,x)$.
\end{thm}

This theorem has two important consequences (Corollaries \ref{cortdos}, \ref{crucial}), that will be crucial for the enumeration of orbits of $n$-sets and $n$-multisets:

\begin{cor}\label{cortdos}
Let $\ga\in \pgl$. Let $W\subseteq V$ be subvarieties of $\pr{N}$ defined over $k$, both expressable as a finite union of linear irreducible $\ga$-invariant subvarieties of $\pr{N}$. Let $U=V\setminus W$ be the complementary variety. Then,    
$$
f_{V/\ga}(x)=f_V(x),\quad f_{U/\ga}(x)=f_U(x),\quad \overline{f}_{V/\ga}(x)=\overline{f}_V(x),\quad \overline{f}_{U/\ga}(x)=\overline{f}_U(x).
$$
\end{cor}
\begin{pf}
If $V$ is a linear irreducible $\ga$-invariant subvariety of $\pr{N}$, then $V\simeq \pr{\,\dm V}$ and $\op{Z}\left((V/\ga)/k,x\right)=\op{Z}(V/k,x)$ by
Theorem \ref{tdos}. This equality holds too for $V$ a finite union of linear irreducible $\ga$-invariant subvarieties, since each irreducible component of $V$ and the intersection of an arbitrary number of components are projective spaces. The corollary follows then from Theorem \ref{tu} and Corollary \ref{multipl}.\hfill{$\Box$}
\end{pf}

\begin{defn}
For each $P\in\pr{N}$ we denote by $O_{\ga}(P)$ the orbit of $P$ under the action of the cyclic group generated by $\ga$.  This set will be simply called  the ``$\ga$-orbit of $P$".

Let $\ga\in\pgl$ and let $V\subseteq \pr{N}$ be a subvariety defined over $k$. We shall use the following notation for the respective sets of rational $n$-sets and $n$-multisets of $V$ that are fixed by $\ga$ as unordered families of points of $\pr{N}(\kb)$:
\begin{equation}\label{fixgamma}\as{2.2}
\begin{array}{l}
\fix{V}{n}:=\{S\in \comb{V}{n}(k)\tq \ga(S)=S\},\\ \fixr{V}{n}:=\{S\in \combr{V}{n}(k)\tq \ga(S)=S\}.
\end{array}
\end{equation}
\end{defn}

\begin{cor}\label{crucial}
Let $\ga$, $V$, $W$, $U$ be as in Corollary \ref{cortdos} and suppose that all $\kb$-rational points of $U$ have $\ga$-orbits of the same length $m$. Then,
$$
|\fix{U}{mn}|=\left|\comb{U}{n}(k)\right|,\qquad 
|\fixr{U}{mn}|=\left|\combr{U}{n}(k)\right|.
$$
\end{cor}
\begin{pf}
The $\ga$-invariant and $\sigma$-invariant $mn$-sets of $U$ are in one-to-one 
correspondence with the $\sigma$-invariant $n$-sets of $U/\ga$; thus,  by Corollary \ref{cortdos}
$$|\fix{U}{mn}|=\left|\comb{U/\ga}{n}(k)\right|=\left|\comb{U}{n}(k)\right|.$$ The argument for 
$|\fixr{U}{mn}|$ is analogous.\hfill{$\Box$}\end{pf}
In order to prove Theorem \ref{tdos} we show first a similar result for the affine space.
\begin{prop}\label{pdotze}
For any $\ga\in\gl{N}{k}$ we have $|(\af{N}/\ga)(k)\vert=q^N$.
\end{prop}
\begin{pf}
Our aim is to compute the cardinality of the set 
$$
(\af{N}/\ga)(k)=\{O_\ga(P), P\in\af{N}(\kb)\tq  \sigma(O_\ga(P))=O_\ga(P)\}.
$$
Since $\ga$ and $\sg$ commute, for any $P\in\af{N}(\kb)$ we have:
\begin{equation}\label{lem1}
\sigma(O_\ga(P))=O_\ga(P)\ \mbox{ if and only if }\ \sigma(P) \in O_\ga(P).
\end{equation}
For any $\rho\in \gl{N}{k}$, let us denote by $C_\rho$ the set $\{P\in\af{N}(\kb)\tq  \sigma(P)=\rho(P)\}$. If $m$ is the order of $\ga$ as an element of $\gl{N}{k}$ we claim that
\begin{equation}\label{lem2}
|(\af{N}/\ga)(k)|=\frac1m\sum_{0\leq i<m}|C_{\ga^i}|.
\end{equation}
In fact, consider the formal disjoint union of all $C_{\gamma^i}$ (they are not disjoint as subsets of $\af{N}(\kb)$) and the map
$$
\renewcommand\arraystretch{1.4}
\begin{array}{rccl}
O_\ga\colon\negmedspace& \coprod_{0\le i<m}C_{\ga^i}&\lra&(\af{N}/\ga)(k)\\
&P&\mapsto&O_{\ga}(P)
\end{array}
\renewcommand\arraystretch{1.}
$$
By (\ref{lem1}), $O_\gamma(P)$ is defined over $k$ if and only if
 $P\in\cup_{0\leq i<m}C_{\gamma^i}$, so that this map is well-defined and onto. Thus, to prove (\ref{lem2}) we need only to check that
  each $\ga$-orbit $O_\gamma(P)\in(\af{N}/\gamma)(k)$ has  exactly $m$ preimages. Let $D=|O_\ga(P)|$; clearly $D|m$ and from $\sigma(P)\in
O_\ga(P)=\{P,\ga(P),\dots,\ga^{D-1}(P)\}$, we see that $P \in C_{\gamma^i}$ for a unique $0\le i<D$. On the other hand, $O_\ga(P)\subseteq C_{\ga^i}$
because
$$
\sg(\ga^j(P))=\ga^j(\sg(P))=\ga^{j+i}(P)=\ga^i(\ga^j(P)).
$$
Hence, the $D$ elements in $O_\gamma(P)$ are precisely the preimages of $O_\gamma(P)$ by the map $O_\gamma$ restricted to $C_{\gamma^i}$. Now, all these
points belong to 
 $$
C_{\gamma^{i+D}},C_{\gamma{i+2D}},\dots,C_{\gamma^{i+(\frac
mD-1)D}},
 $$and none of these points belongs to any other $C_{\gamma^j}$. Therefore, $O_\gamma(P)$ has
 exactly $D(m/D)=m$ preimages.

Finally, the proposition will be proved if we show that $|C_\rho|=q^N$ for all $\rho\in \gl{N}{k}$. Let us check this; for any given $\rho\in \gl{N}{k}$ let $\beta\in\gl{N}{k}$ be such that $\beta\rho\beta^{-1}$ is a rational canonical matrix:  $\beta\rho\beta^{-1}=\di(A_1,\dots,A_r)$, each $A_i$ being a cyclic component of the type:
\begin{equation}\label{edotze}
 \as{1}A=\begin{pmatrix}
0&&&&-a_s\\1&0&&&-a_{s-1}\\&1&\ddots&&\vdots\\&&\ddots&0&-a_2\\&&&1&-a_1
 \end{pmatrix}, 
 \end{equation}
with $x^s+a_1x^{s-1}+\cdots+a_s\in k[x]$ an invariant factor of the endomorphism $\rho$. The sets $C_{\rho}$ and $C_{\beta\rho\beta^{-1}}$ have the same cardinality because the automorphism $\beta$ of $\af{N}(\kb)$ maps one set onto the other. Now, if we split the coordinates of the points $x\in\af{N}(\kb)$ into $x=(x_1,x_2,\dots,x_r)$ with $x_i\in \kb^{s_i}$, $s_1+\cdots+s_r=N$,
the condition $\sigma(x)=\beta\rho\beta^{-1}(x)$ translates into
  $\sigma(x_i)= A_ix_i$ for all $i=1,\dots,r$; thus, we need only to check that $|C_A|=q^s$, for $A$ a companion matrix as in (\ref{edotze}).
For $x\in \kb^s$, the equality $\sigma(x)=A(x)$ splits into
\begin{equation}\label{equinze}
\sigma(x_1)=-a_sx_s,\quad \sigma(x_2)=x_1-a_{s-1}x_s,\  \cdots\ \sigma(x_s)=x_{s-1}-a_1x_s.
\end{equation}
This allows us to express $x_1,\dots,x_{s-1}$ as a linear combination of $x_s$ and its Galois conjugates:
$$
x_{s-i}=a_ix_s+a_{i-1}\sigma(x_s)+\cdots+a_1\sigma^{i-1}(x_s)+\sigma^i(x_s),\quad 1\le i<s,
$$
hence, the first equation of (\ref{equinze}) is equivalent to
$$
\sigma^s(x_s)+a_1\sigma^{s-1}(x_s)+\cdots+a_{s-1}\sigma(x_s)+a_sx_s=0.
$$
Since $a_s\ne0$, this is a separable equation in $x_s$ with $q^s$ solutions in $\kb$.\hfill{$\Box$}
\end{pf}
We are now ready to prove Theorem \ref{tdos}. Actually it is an immediate consequence of the following result:
\begin{prop}
For any $\ga\in\pgl$ we have $|(\pr{N}/\ga)(k)|=\dfrac{q^{N+1}-1}{q-1}$.
\end{prop}
\begin{pf}
We choose a representative of $\ga$ in $\gl{N+1}{k}$, which we still denote by
 $\ga$. In the sequel we identify an affine point $P\in\af{N+1}(\kb)$ with its image $P\in\pr{N}(\kb)$ under 
 the  natural morphism $\pi\colon \af{N+1}\setminus\{0\}\lra
\pr{N}$. However, in order to avoid confusion we shall denote by $O_\ga(P)$ the affine $\ga$-orbit of $P$ and by $O_{\ga}^{\mbox{\scriptsize $\op{pr}$}}(P)$ the projective orbit. Clearly $\pi(O_\ga(P))=O_{\ga}^{\mbox{\scriptsize $\op{pr}$}}(P)$ and $\pi$ induces a natural map
$$
\pi: (\af{N+1}/\ga)(k)\setminus\{0\}\lra (\pr{N}/\ga)(k).
$$By Proposition \ref{pdotze}, in order to prove the proposition we need only to show that $\pi$ is onto and each element of $(\pr{N}/\ga)(k)$ has  $q-1$ preimages in $(\af{N+1}/\ga)(k)$.

For $P,\,Q\in\af{N+1}(\kb)\setminus\{0\}$, the condition $\pi(O_\ga(Q))=O_\ga^{\mbox{\scriptsize $\op{pr}$}}(P)$ is equivalent to
$$
\exists i\in\N,\,
\mu\in\kb^*:\ Q=\mu\gamma^i(P), \ \mbox{ or equivalently }\ \exists\mu\in\kb^*:\
O_\gamma(Q)=O_\gamma(\mu P).
$$
Thus, if for some $P$ the orbit $O_{\gamma}^{\mbox{\scriptsize $\op{pr}$}}(P)$ is defined over 
$k$ we want to check that exactly $q-1$ of the orbits
$O_\gamma(\mu P)$, $\mu\in\kb^*$,  are defined over $k$. To check this, consider the following subgroup 
of $\kb^*$:
$$
\Lambda_P:=\{\lambda\in\kb^*\tq\lambda O_\gamma(P)=O_\gamma(P)\},
$$and let $e=|\Lambda_P|$. For any $\mu\in\kb^*$,
$$
\sigma(O_\gamma(\mu P))=\sigma(\mu
O_\gamma(P))=\sigma(\mu)\sigma(O_\gamma(P))=\sigma(\mu)O_\gamma(P),
$$
and this coincides with $O_\gamma(\mu P)=\mu O_\gamma(P)$ if and only if
 $\sigma(\mu)\mu^{-1}O_\gamma(P)=O_\gamma(P)$. Hence, the orbit $O_\gamma(\mu P)$ is defined over $k$
precisely for the $e(q-1)$ values of $\mu\in\kb^*$ determined by the condition $\mu^{q-1}\in \Lambda_P$.
For any such $\mu$  there are $e$ other values with the same affine $\ga$-orbit, because
$$
 O_\gamma(\mu P)= O_\gamma(\mu' P)
 \Leftrightarrow\mu'\mu^{-1}\in\Lambda_P\Leftrightarrow\mu'\in\mu\Lambda_P,
$$
Therefore, among all $O_\gamma(\mu P)$ there are 
$q-1$ different $\ga$-orbits that are defined over $k$.\hfill{$\Box$}
\end{pf}
\section{Proper subvarieties of $\pr{N}$ with respect to a fixed automorphism}\label{sectres}
Throughout this section we fix a $k$-automorphism of $\pr{N}$, represented by some $\ga\in\pgl$.

Let $\ll_{\ga}$ be the poset of $\ga$-invariant irreducible linear subvarieties of $\pr{N}(\kb)$, ordered by inclusion.
The poset $\ll_{\ga}$ is not locally finite. For instance, if $V$ is a plane of fixed points of $\ga$ and $P\in V$,
the interval $[P,V]$ is not finite.
For the basic concepts and notations about posets we address the reader to \cite{sta}.
\begin{defn}
(1) For any $V\in\ll_{\ga}$ we define the exponent of $V$ as the order of $\ga$ as a projective automorphism of $V$ 
$$
\op{exp}V:=\op{ord}(\gamma_{|V}).
$$Note that for all $V,\,W\in\ll_{\ga}$: $\ V\le W\imp \exp V\,|\,\exp W$.

(2) A node $V\in\ll_{\ga}$ is said to be proper if it is maximal among all nodes with the same exponent:
$$\exp V<\exp W,\quad\forall W\in \ll_{\ga}\mbox{ such that } V<W.$$
We denote by $\ig$ the subposet of $\ll_{\ga}$ formed by the proper nodes.
\end{defn}
Note that $\pr{N}$ is always proper. In our example above, the points $P$ and lines $L$ of a plane $V$ of fixed points of $\ga$ are not proper, because they have the same exponent than $V$: $\exp P=\exp L=\exp V=1$. 

In this section we shall see that $\ig$ is a finite poset and we shall determine its structure. 
To this end we need to introduce some terminology. We fix a representative of $\ga$ in $\gl{N+1}{k}$, which we still denote by $\ga$; we abuse of language and use the same notation for $\ga$-invariant subvarieties of $\pr{N}(\kb)$ and their affine cones, which are $\ga$-invariant linear subspaces of $\af{N+1}(\kb)$.

Let $\op{VP}_{\ga}=\{\la_1,\dots,\la_s\}$ be the set of eigenvalues of $\ga$. Recall the decomposition
$$
\af{N+1}(\kb)=V_1\oplus\cdots\oplus V_s,\quad  
V_i=\op{Ker}(\gamma-\lambda_i)^{m_i},\ i=1,\dots,s,
$$
where $m_i$ is the maximum exponent such that $(x-\la_i)^{m_i}$ divides the minimal polynomial of $\ga$.
Each $V_i$ is $\ga$-invariant and $\exp V_i=p^{\delta_i}$, with $\delta_i=\lceil\log_p(m_i)\rceil$. 
 
Let $\mathbf{D}_{\N}$ be the poset of positive integers ordered by divisibility. Let $\bbb_{\ga}$ be the poset of nonempty subsets of $\op{VP}_{\ga}$
ordered by inclusion. For each $\Lambda\in\bbb$ we define two invariants, $\delta(\Lambda)$, $D(\Lambda)$ in the form of two morphisms of posets
\begin{equation}\label{deltad}
\delta\colon\bbb_{\ga}\lra\N,\qquad D\colon \bbb_{\ga}\lra\mathbf{D}_{\N}, 
\end{equation}
$$\delta(\Lambda):=\max\{\lceil\log_p(m_i)\rceil\tq \la_i\in\Lambda\}=\max\{\delta_i\tq \la_i\in\Lambda\},$$
$$D(\{\lambda_{i_1},\dots,\lambda_{i_t}\}):=\ord(\di(\lambda_{i_1},\dots,\lambda_{i_t}))=\op{min}\{d
\tq\lambda_{i_1}^d=\lambda_{i_2}^d=\cdots=\lambda_{i_t}^d\}.$$
Note that $D(\Lambda)$ is always prime to $p$. 
\begin{defn}
We say that $\Lambda\in\bbb_{\ga}$ is $D$-proper if $\Lambda$ is maximal among all nodes with the same value of $D$:
$$D(\Lambda)<D(\Delta),\quad\forall \Delta\in \bbb_{\ga}\mbox{ such that } \Lambda<\Delta.$$
We denote by $\bbb_{\ga}^{\op{pr}}$ the subposet of the $D$-proper nodes of $\bbb_{\ga}$. 
\end{defn}
The following remark is obvious:
\begin{lem}
The $\gamma$-invariant linear spaces of  $\af{N+1}(\kb)$ are all of the form $W=W_1\oplus\cdots\oplus W_s$,
with $W_i\subseteq V_i$ $\gamma$-invariant.  If each $W_i$ has exponent $\exp W_i=p^{\eps_i}$, $\eps_i\le \delta_i$, then $\exp W=p^{\eps}D(\Lambda_W)$, where $\eps=\max\{\eps_i\}$ and $\Lambda_W:=\{\la_i\in\op{VP}_{\ga}\tq W_i\ne0\}$.
\end{lem}
For any $0\le \nu\le\delta_i$ the subspace $\k(\ga-\la_i)^{p^{\nu}}$ is the maximum $\ga$-invariant subspace of $V_i$ with exponent $p^{\nu}$. The following result follows immediately:
\begin{lem}\label{maxim}
For each pair $(\Lambda,\nu)$ with $\Lambda\in\bbb_{\ga}$ and $0\le\nu\le \delta(\Lambda)$, consider the $\ga$-invariant linear subvariety of $\pr{N}(\kb)$:  
$$V_{\Lambda}^{(\nu)}:=\bigoplus_{\la\in\Lambda}\op{Ker}(\ga-\la)^{p^{\nu}}.$$
Then $\exp(V_{\Lambda}^{(\nu)})=p^{\nu}D(\Lambda)$, and this subvariety $V_{\Lambda}^{(\nu)}$ contains all $\ga$-invariant subvarieties $W$ such that $\Lambda_W\subseteq\Lambda$ and $v_p(\exp W)\le \nu$.  
\end{lem}
Note that for any $\Lambda\in\bbb_{\ga}$ with invariants $D=D(\Lambda)$, $\delta=\delta(\Lambda)$, we have a chain
of nodes of the poset $\ll_{\ga}$ with respective exponents $D,\,pD, \dots, p^{\delta}D$:
$$
V_{\Lambda}^{(0)}<V_{\Lambda}^{(1)}<\cdots<V_{\Lambda}^{(\delta)}.
$$ 
\begin{thm}\label{pvint}
A node $V\in\ll_{\ga}$ is proper if and only if $V=V_{\Lambda}^{(\nu)}$ for some $\Lambda\in \bbb_{\ga}$ which is $D$-proper, and some $0\le\nu\le\delta(\Lambda)$.
\end{thm}
\begin{pf}
Suppose $\Lambda$ is $D$-proper and let us show that $V_{\Lambda}^{(\nu)}$ is proper.  Suppose that $V_{\Lambda}^{(\nu)}\subseteq W$ for some $W\in\ll_{\ga}$ with $\exp(V_{\Lambda}^{(\nu)})=\exp(W)$. 
In particular $D(\Lambda)=D(\Lambda_W)$ and this implies $\Lambda=\Lambda_W$ because $\Lambda$ is $D$-proper. Hence, $W=V_{\Lambda}^{(\nu)}$ by Lemma \ref{maxim}. 

Conversely, suppose  $W$ proper with $\exp W=p^{\eps}D(\Lambda_W)$. By Lemma \ref{maxim} we have $W\subseteq V_{\Lambda_W}^{(\eps)}$ and this implies $W= V_{\Lambda_W}^{(\eps)}$ because $W$ is proper. Finally, $\Lambda_W$ is proper because $\Lambda_W\subsetneq\Lambda$, with $D(\Lambda_W)=D(\Lambda)$ would lead to $W= V_{\Lambda_{W}}^{(\eps)}\subsetneq V_{\Lambda}^{(\eps)}$ and $W$ would not be  proper. \hfill{$\Box$}
\end{pf}
\begin{cor}\label{pvintcor}
Let $\delta=\delta(\op{VP}_{\ga})$ and let $[0,\delta]$ be the poset of integers $0\le \nu\le \delta$ ordered by size. We have a natural identification: 
$$
 \ig\hookrightarrow [0,\delta]\,\times\,\bbb_{\ga}^{\op{pr}},
$$
with image the subposet containing the nodes $(\nu,\Lambda)$ with $\nu\le\delta(\Lambda)$. In parti\-cular, $\ig$ is a finite poset.
Moreover, $\ig$ is a lattice.
\end{cor}
\begin{pf}
The first statement is an immediate consequence of Theorem \ref{pvint}. In order to prove that $\ig$ is a lattice it is sufficient to check that the subposet $\bbb_{\ga}^{\op{pr}}$ is a
lattice. The total set $\hat{1}$ is always proper; hence, by \cite[3.3.1]{sta} it is sufficient to check that the intersection of two 
$D$-proper elements is $D$-proper. Consider pairwise disjoint sets $\Lambda_0,\,\Lambda_1,\,\Lambda_2\in\bbb_{\ga}$, such that  
$\Lambda_0\cup\Lambda_1$ and $\Lambda_0\cup\Lambda_2$ are $D$-proper. Let us show that $\Lambda_0$ is $D$-proper. If  $D(\Lambda_0\cup\{\la\})=D(\Lambda_0)$, then $$D(\Lambda_0\cup\Lambda_1\cup\{\la\})=D(\Lambda_0\cup\Lambda_1),\qquad
D(\Lambda_0\cup\Lambda_2\cup\{\la\})=D(\Lambda_0\cup\Lambda_2).$$Since these sets are $D$-proper we have $\la\in  (\Lambda_0\cup\Lambda_1)\cap(\Lambda_0\cup\Lambda_2)=\Lambda_0$. \hfill{$\Box$}
\end{pf}
\begin{lem}\label{disjoint}
Two proper $\ga$-invariant subvarieties of the same exponent are either disjoint or coincident.
\end{lem}
\begin{pf}Suppose that $V,\,W\in\ig$ have both exponent $d$, and $V\cap W\ne\emptyset$. Then, the linear subvariety generated by $V$ and $W$ is $\ga$-invariant and it has still exponent $d$; by the properness of $V$ and $W$ we have necessarily $V=W$.\hfill{$\Box$}
\end{pf}
We finish this section with a crucial property of the proper $\ga$-invariant subvarieties. Recall that the degree $\dg V$ of a subvariety $V\subseteq \pr{N}(\kb)$ is the minimum exponent $r$ such that $\sg^r(V)=V$, or, equivalently, the minimum positive integer $r$ such that $V$ is defined over $k_r$. 
\begin{prop}\label{zetavg}
Let $V$ be a proper $\ga$-invariant subvariety of $\pr{N}(\kb)$ of degree $r$. Consider the $\ga$-invariant linear variety defined over $k$ $$V_G:=V\cup \sg(V)\cup\cdots\cup\sg^{r-1}(V).$$ Then, $\op{Z}(V_G/k,x)=\op{Z}(V/k_r,x^r)$.   
\end{prop}
\begin{pf}
 We want to prove the identity of formal series:
 \begin{equation}\label{zig}
 \sum_{m\ge1}\frac{|V_G(k_m)|}{m}x^m= \sum_{n\ge1}\frac{|V(k_{rn})|}{n}x^{rn}.
\end{equation}
By Lemma \ref{disjoint} the varieties $\sg^i(V)$ are pairwise disjoint, because they are all proper and have the same exponent. Since,     
$$
P\in V_G(k_m) \ \imp\ P\in \sg^i(V),\ \sg^m(P)=P\ \imp\ P\in \sg^i(V)\cap \sg^{i+m}(V),
$$
we deduce that $V_G(k_m)=\emptyset$ if $r\nmid m$; hence, we can change $m=rn$ in the left side of (\ref{zig}), and the equality holds because $|V_G(k_{rn})|=r|V(k_{rn})|$.\hfill{$\Box$}
\end{pf}

\section{G-exponent index and generating functions}\label{secquatre}
Let $\g$ be a finite group acting on a finite set $X$. The number of orbits of this action can be counted as the average number of fixed points (\cite{bfkwz}, 3.1.6):
\begin{equation}\label{cf} |\g\backslash X|=\frac
1{\vert\g\vert}\sum_{\ga\in\g}\vert \xg\vert=
\sum_{\ga\in\cc}\frac{\vert \xg\vert}{\vert \g_{\ga}\vert},
\end{equation} where  $\cc$ is a set of representatives of conjugacy classes of elements of
$\g$ and
 $$
 \xg:=\{x\in X \tq \ga(x)=x\}, \quad
 \g_{\ga}:=\{\rho\in \g \tq \rho\ga\rho^{-1}=\ga\}.
 $$

In this section we apply this formula to compute the generating function of 
$$t_N(n):=|\g\backslash X|, \ \overline{t}_N(n):=|\g\backslash Y|,\quad \mbox{ for }X=\comb{\pr{N}}{n}(k), \ Y=\combr{\pr{N}}{n}(k),
$$and $\g:=\pgl$. With the notation introduced in (\ref{fixgamma}), the crucial step is the computation of $|\fix{\pr{N}}{n}|$, $|\fixr{\pr{N}}{n}|$. The main ingredient is the quotient poset of $\ig$ under the  Galois action.
 \begin{defn}
For any $\ga\in\pgl$ we denote by $\igg$ the quotient poset of $\ig$ under the  Galois action. That is, for any $V\in\ig$ of degree $r$, the whole $\sg$-orbit $V,\,\sg(V),\,\cdots,\,\sg^{r-1}(V)$ determines a single node of $\igg$, denoted by $[V]$. Each node $[V]\in\igg$ determines a $\ga$-invariant linear subvariety $V_G:=V\cup \sg(V)\cup\cdots\cup\sg^{r-1}(V)$, which is defined over $k$ and is independent of the choice of $V$ inside a given $\sg$-orbit. 
The poset $\igg$ inherites the following ordering of $\ig$: 
$$[W]\le[V]\sii W_G\subseteq V_G\sii W\le\sg^i(V),\ \mbox{ for some }i.
$$
Usually we shall abuse of terminology and write simply $V\in\igg$ to indicate the class in $\igg$ of certain $V\in\ig$.  
\end{defn}
For any $V\in\ig$ we define
$$
V^0:=V\setminus\left(\bigcup_{W\in\ig,\,W<V}W\right),\qquad V_G^0:=V_G\setminus\left(\bigcup_{W\in\igg,\,W<V}W_G\right).
$$
Clearly, $V_G^0=V^0\cup \sg(V^0)\cup\cdots\cup\sg^{r-1}(V^0)$, where $r=\dg V$.
\begin{lem}\label{gaorbit}
For each $P\in V^0$, the $\ga$-orbit of $P$ has $\exp V$ elements. 
\end{lem}
\begin{pf}
The length $|O_{\ga}(P)|$ of the $\ga$-orbit of any $P\in V^0$ is the minimum positive integer $e$
such that $P$ is a fixed point of $\ga^e$. Now, the linear subvariety $W$ generated by $O_{\ga}(P)$ is 
 $\ga$-invariant and it is pointwise fixed by $\ga^e$. Hence, $e=\ord(\ga_{|W})=\exp W$.

Embed $W$ in a (unique) proper subvariety with the same exponent:
 $$
 W\subseteq W',\quad \exp W=\exp W',\quad W'\in \ig.
 $$ 
Since $V$ is proper and $\ig$ is a lattice, we have necessarily $W\subseteq W'\subseteq V$. Finally, since $P\in V^0$, the point $P$ cannot lie in a proper strict subvariety of $V$; thus, $W'=V$ and $e=\exp W=\exp V$. \hfill{$\Box$}
\end{pf}
Consider the stratification $\pr{N}(\kb)=\coprod_{V\in \igg}V_G^0$. For any $n$-set $S\in \fix{\pr{N}}{n}$ let the distribution of the $n$ points of $S$ among these strata be 
$$
S=\coprod_{V\in \igg}S_V,\quad S_V=S\cap V_G^0.
$$Since $S$ and $V_G^0$ are $\ga$-invariant and $\sg$-invariant, each $S_V$ is a $\ga$-invariant and $\sg$-invariant unordered family of points; in other words, $S_V\in\fix{V_G^0}{|S_V|}$. Thus, we can count the number of possibilities for $S$ just by considering all possible numerical distributions of $n$ points among the strata $V_G^0$, and then counting, for each numerical distribution, the number of possibilities for $S_V$. By Lemma \ref{gaorbit},  $|S_V|=n_V\exp V$ for some nonnegative integer $n_V$, and we get
\begin{multline*}
|\fix{\pr{N}}{n}|=\sum_{\sum_{V\in
\igg} n_V\op{exp} V=n}\left(\prod_{V\in
\igg}\left|\fix{V_G^0}{n_V\exp V}\right|\right)=\\=\sum_{\sum_{V\in
\igg} n_V\op{exp} V=n}\left(\prod_{V\in
\igg}\left|\comb{V_G^0}{n_V}(k)\right|\right),
\end{multline*}
the last equality by Lemma \ref{gaorbit} and Corollary \ref{crucial}.
Therefore, the generating function of these numbers is: 
$$
\sum_{n\ge0}|\fix{\pr{N}}{n}| x^n=\prod_{V\in\igg}f_{V_G^0}(x^{\exp V})=\prod_{V\in\igg}f_{V^0}(x^{\exp V\dg V}),
$$the last equality by Proposition \ref{zetavg} and Theorem \ref{tu}.
An application of the Cauchy-Frobenius formula (\ref{cf}) leads to a first computation of the generating function we are interested in:
\begin{equation}\label{first}\as{1.8}
\sum_{n\ge 0}t_N(n)x^n=\sum_{\ga\in\cc}|\g_{\ga}|^{-1}\prod_{V\in\igg}f_{V^0}(x^{\exp V\dg V}).
\end{equation}
We can refine (\ref{first}) by grouping together all elements $\ga\in\cc$ with a common value of  $\prod_{V\in\igg}f_{V^0}(x^{\exp V\dg V})$.
By Moebius inversion in the poset $\igg$ 
$$f_{V^0}(x)=\prod_{W\le V}f_W(x)^{\mu(W,V)},$$
and the function $f_W(x)$ depens only on $\dm W$, because $W\simeq\pr{\,\dm W}$. Hence, the term $\prod_{V\in\igg}f_{V^0}(x^{\exp V\dg V})$
depens only on the structure of the poset $\igg$ and the triple weight $(\dm V ,\exp V ,\dg V)$ of each node. 
\begin{defn}
We say that two elements $\ga,\ga'\in\cc$  have the same {\it subtype}, and we write $\ga\sim \ga'$, if there exists a poset isomorphism $\igg\,\iso\, \ll_{\ga',G}^{\op{pr}}$, preserving the weight $(\dim V,\exp V,\dg V)$ of each node. We denote by $\mathcal{S}:=\cc/\sim$ the quotient set of $\cc$ by this equivalence relation, and by $\st\colon \cc\to\mathcal{S}$ the canonical quotient map. 
For each subtype $\al\in\mathcal{S}$ we denote by $\ll(\al)$ the poset $\igg$ for any $\ga$ with $\st(\ga)=\al$, and we consider the weighted sum:
$$
M_{\al}:=\sum_{\st(\ga)=\al}|\g_{\ga}|^{-1}.
$$
\end{defn}
Our main theorem is a rewriting of (\ref{first}) after grouping together all $\ga\in\cc$ in the same subtype. We include in the theorem the similar statement for $n$-multisets, which is obtained by completely analogous arguments.
\begin{thm}\label{ffinal}
$$\sum_{n\ge0}t_N(n)x^n=\sum_{\al\in\mathcal{S}}M_{\al}\prod_{V\in
\ll_G(\al)}f_{V^0}(x^{\exp(V)\op{deg}(V)}).
$$
$$\sum_{n\ge0}\overline{t}_N(n)x^n=\sum_{\al\in\mathcal{S}}M_{\al}\prod_{V\in
\ll_G(\al)}\overline{f}_{V^0}(x^{\exp(V)\op{deg}(V)}).
$$
\end{thm}

This formula is suitable of an effective implementation. In this regard one needs only to carry out the following tasks:
\begin{enumerate}
\item Find an intrinsic description of the set $\mathcal{S}$.
\item For each $\al\in\mathcal{S}$ find an intrinsic description of the weighted poset $\ll(\al)$ and its Moebius function.
\item For each $\al\in\mathcal{S}$ find an explicit formula for the universal coefficients $M_{\al}$.
\end{enumerate}
This will be fulfilled in section \ref{seccinc} for the cases $N=1,\,2$. As a consequence, one is able to deduce explicit formulas for $t_2(n)$, $\overline{t}_2(n)$ as polynomials with integer coefficients in the cardinality $q$ of the ground field. Similar formulas for $t_1(n), \,\overline{t}_1(n)$ had been obtained in \cite{lmnx}. 

\section{Explicit formulas for dimension $N=1,\,2$}\label{seccinc}
\begin{scriptsize}
\begin{center}
\begin{longtable}{|c|c|c|c|}
\hline
type&$\ga$&$\ll(\al)=\igg$&subtypes\\
\hline
A&$B_2$&\setlength{\unitlength}{4.mm}
\lower3.ex\hbox{\begin{picture}(2,2.4)
\put(1,1.3){$\bullet$}\put(1,0.3){$\circ$}
\put(1.5,1.4){\begin{scriptsize}$d$\end{scriptsize}}
\put(1.5,.4){\begin{scriptsize}$1\,(2)$\end{scriptsize}}
\put(1.15,.62){\line(0,1){1}}
\end{picture}}&$d\tq q+1,\ d>1$\\
\hline B&\as{1}
$\begin{pmatrix}1&0\\1&1\end{pmatrix}$&
\setlength{\unitlength}{4.mm}
\lower1.ex\hbox{\begin{picture}(2,2)
\put(1,0){$\bullet$}\put(1,1){$\bullet$}
\put(1.5,.1){\begin{scriptsize}$1$\end{scriptsize}}
\put(1.5,1.1){\begin{scriptsize}$p$\end{scriptsize}}
\put(1.2,.2){\line(0,1){1}}
\end{picture}}&\\\hline C&
$\as{.8}\begin{array}{c}\di(1,\la)\\\la\ne1\end{array}$&
\setlength{\unitlength}{4.mm}
\lower2.ex\hbox{\begin{picture}(2,2.4)
\put(0,0.3){$\bullet$}\put(1,1.3){$\bullet$}\put(1,.3){$\bullet$}
\put(-.4,.4){\begin{scriptsize}$1$\end{scriptsize}}
\put(1.5,.4){\begin{scriptsize}$1$\end{scriptsize}}
\put(1.5,1.4){\begin{scriptsize}$d$\end{scriptsize}}
\put(1.2,.5){\line(0,1){1}}
\put(.2,.5){\line(1,1){1}}
\end{picture}}&$d\tq q-1,\ d>1$
\\\hline D&
$1$&\setlength{\unitlength}{4.mm}
\begin{picture}(1,1)
\put(0,0){$\bullet$}
\put(.6,.1){\begin{scriptsize}$1$\end{scriptsize}}
\end{picture}&\\\hline
\caption{Types, subtypes and posets $\ll(\al)$ for $\ga\in\op{PGL}_2(k)$}
\end{longtable}         
\end{center}
\begin{longtable}{|c|c|c|c|}
\hline
type&$\ga$&$\ll(\al)=\igg$&subtypes\\
\hline
A&$B_3$&\setlength{\unitlength}{5.mm}
\begin{picture}(2,1.6)
\put(1,.9){$\bullet$}\put(1,-.6){$\circ$}
\put(-1.6,-.65){\begin{scriptsize}$\dm\,0$\end{scriptsize}}
\put(1.5,.9){\begin{scriptsize}$d$\end{scriptsize}}
\put(1.5,-.6){\begin{scriptsize}$1\,(3)$\end{scriptsize}}
\put(1.15,-.4){\line(0,1){1.5}}
\end{picture}&$\begin{array}{c}d\tq q^2+q+1\\ d>1\end{array}$\\
\hline B&
$\begin{array}{c}\\\di(B_2,1)\\\\\end{array}$&
\setlength{\unitlength}{5.mm}
\begin{picture}(5,2.2)
\put(-1,-.6){$\bullet$}
\put(0,1.4){$\bullet$}\put(0,-.6){$\circ$}\put(0,.4){$\bullet$}
\put(-1.4,-.5){\begin{scriptsize}$1$\end{scriptsize}}
\put(.5,.5){\begin{scriptsize}$d$\end{scriptsize}}
\put(.5,1.5){\begin{scriptsize}$de$\end{scriptsize}}
\put(.5,-.5){\begin{scriptsize}$1\,(2)$\end{scriptsize}}
\put(.17,-.3){\line(0,1){1.9}}\put(-.8,-.4){\line(1,2){.9}}
\put(-1,-1.2){$e>1$}
\put(4,-.6){$\bullet$}
\put(5,1.4){$\bullet$}\put(5,-.6){$\circ$}
\put(3.6,-.5){\begin{scriptsize}$1$\end{scriptsize}}
\put(5.5,1.5){\begin{scriptsize}$d$\end{scriptsize}}
\put(5.5,-.5){\begin{scriptsize}$1\,(2)$\end{scriptsize}}
\put(5.15,-.3){\line(0,1){1.9}}\put(4.2,-.3){\line(1,2){.9}}
\put(4,-1.2){$e=1$}
\end{picture}&$\begin{array}{l}
d\tq q+1,\ d>1\\e\tq q-1
\end{array}$\\\hline \raise4ex\hbox{C}&\as{1}
\raise4.5ex\hbox{$\as{.9}\begin{pmatrix}1&0&0\\1&1&0\\0&1&1\end{pmatrix}$}&
\setlength{\unitlength}{5.mm}
\lower1.ex\hbox{\begin{picture}(3.2,4)
\put(0,1){$\bullet$}\put(0,3){$\bullet$}
\put(3,1){$\bullet$}\put(3,2){$\bullet$}\put(3,3){$\bullet$}
\put(.5,1.1){\begin{scriptsize}$1$\end{scriptsize}}
\put(.5,3.1){\begin{scriptsize}$p$\end{scriptsize}}
\put(3.5,1.1){\begin{scriptsize}$1$\end{scriptsize}}
\put(3.5,2.1){\begin{scriptsize}$2$\end{scriptsize}}
\put(3.5,3.1){\begin{scriptsize}$4$\end{scriptsize}}
\put(-.4,.3){\begin{scriptsize}$p>2$\end{scriptsize}}
\put(2.6,.3){\begin{scriptsize}$p=2$\end{scriptsize}}
\put(0.17,1.2){\line(0,1){2}}
\put(3.165,1.2){\line(0,1){2}}
\end{picture}}&\\\hline D&\as{1}
$\di\left(\begin{pmatrix}1&0\\1&1\end{pmatrix},1\right)$&
\setlength{\unitlength}{5.mm}
\lower4.ex\hbox{\begin{picture}(2,3)
\put(1,1){$\bullet$}\put(1,2){$\bullet$}
\put(-1,1.1){\begin{scriptsize}$\dm\,1$\end{scriptsize}}
\put(1.5,1.1){\begin{scriptsize}$1$\end{scriptsize}}
\put(1.5,2.1){\begin{scriptsize}$p$\end{scriptsize}}
\put(1.165,1.2){\line(0,1){1}}
\end{picture}}&\\\hline
E&\as{1}
$\begin{array}{c}\di\left(\!\begin{pmatrix}1&0\\1&1\end{pmatrix},\la\right)\\\la\ne1\end{array}$&
\setlength{\unitlength}{5.mm}
\lower3.ex\hbox{\begin{picture}(2,3.4)
\put(0,.3){$\bullet$}\put(0,1.3){$\bullet$}
\put(1,.3){$\bullet$}\put(1,1.3){$\bullet$}\put(1,2.3){$\bullet$}
\put(-.5,0.4){\begin{scriptsize}$1$\end{scriptsize}}
\put(1.5,.4){\begin{scriptsize}$1$\end{scriptsize}}
\put(-.5,1.4){\begin{scriptsize}$p$\end{scriptsize}}
\put(1.5,1.4){\begin{scriptsize}$d$\end{scriptsize}}
\put(1.5,2.4){\begin{scriptsize}$pd$\end{scriptsize}}
\put(0.17,.5){\line(0,1){1}}
\put(1.17,.5){\line(0,1){2}}
\put(0.2,.5){\line(1,1){1}}
\put(0.2,1.5){\line(1,1){1}}
\end{picture}}&$d\tq q-1,\ d>1$
\\\hline F&
$\as{.8}\begin{array}{c}\di(1,1,\la)\\\la\ne1\end{array}$&
\setlength{\unitlength}{4.mm}
\lower3.ex\hbox{\begin{picture}(2,3.4)
\put(0,1.3){$\bullet$}\put(1,2.3){$\bullet$}\put(1,.3){$\bullet$}
\put(-.4,1.4){\begin{scriptsize}$1$\end{scriptsize}}
\put(1.5,.4){\begin{scriptsize}$1$\end{scriptsize}}
\put(1.5,2.4){\begin{scriptsize}$d$\end{scriptsize}}
\put(1.2,.5){\line(0,1){2}}
\put(.2,1.5){\line(1,1){1}}
\end{picture}}&$d\tq q-1,\ d>1$
\\\hline G&
$\begin{array}{c}\\\\\di(1,\la,\mu)\\\la,\mu\ne1,\,\la\ne\mu\\\\\end{array}$&
\setlength{\unitlength}{6.mm}
\begin{picture}(8.4,2)
\put(0,-.7){$\bullet$}\put(1,-.7){$\bullet$}\put(2,-.7){$\bullet$}
\put(0,.3){$\bullet$}\put(1,.3){$\bullet$}\put(2,.3){$\bullet$}
\put(1,1.3){$\bullet$}
\put(0.3,-.8){\begin{scriptsize}$1$\end{scriptsize}}
\put(1.3,-.8){\begin{scriptsize}$1$\end{scriptsize}}
\put(2.3,-.8){\begin{scriptsize}$1$\end{scriptsize}}
\put(-.3,.5){\begin{scriptsize}$d$\end{scriptsize}}
\put(1.,0){\begin{scriptsize}$e$\end{scriptsize}}
\put(2.3,.4){\begin{scriptsize}$f$\end{scriptsize}}
\put(1.4,1.4){\begin{scriptsize}$m$\end{scriptsize}}
\put(.15,-.5){\line(0,1){1}}\put(2.15,-.5){\line(0,1){1}}\put(1.15,.5){\line(0,1){1}}
\put(1.08,-.57){\line(-1,1){.9}}\put(1.2,-.5){\line(1,1){.9}}
\put(1.12,.44){\line(-1,-1){.92}}\put(1.08,.43){\line(1,-1){.92}}
\put(1.12,1.44){\line(-1,-1){.92}}\put(1.08,1.43){\line(1,-1){.92}}
\put(.1,-1.5){$d>e>f$}
\put(3,-.7){$\bullet$}\put(4,-.7){$\bullet$}\put(5,-.7){$\bullet$}
\put(5,.3){$\bullet$}\put(4,1.3){$\bullet$}
\put(3.3,-.8){\begin{scriptsize}$1$\end{scriptsize}}
\put(4.3,-.8){\begin{scriptsize}$1$\end{scriptsize}}
\put(5.3,-.8){\begin{scriptsize}$1$\end{scriptsize}}
\put(5.3,.4){\begin{scriptsize}$f$\end{scriptsize}}
\put(4.4,1.4){\begin{scriptsize}$m$\end{scriptsize}}
\put(5.16,-.5){\line(0,1){1}}
\put(4.2,-.5){\line(1,1){.9}}\put(4.1,1.43){\line(1,-1){.92}}
\put(3.1,-.5){\line(1,2){.92}}
\put(3.1,-1.5){$d=e>f$}
\put(6,-.7){$\bullet$}\put(7,-.7){$\bullet$}\put(8,-.7){$\bullet$}
\put(7,1.3){$\bullet$}
\put(6.3,-.8){\begin{scriptsize}$1$\end{scriptsize}}
\put(7.3,-.8){\begin{scriptsize}$1$\end{scriptsize}}
\put(8.3,-.8){\begin{scriptsize}$1$\end{scriptsize}}
\put(7.5,1.4){\begin{scriptsize}$m$\end{scriptsize}}
\put(8.12,-.5){\line(-1,2){.92}}\put(7.15,-.5){\line(0,1){1.92}}
\put(6.17,-.5){\line(1,2){.92}}
\put(6.1,-1.5){$d=e=f$}
\end{picture}&$\begin{array}{c}\\\\(d,\,e,\,f)\in\mathcal{S}_G\\
m=\lcm(d,e)\\\\\end{array}$
\\\hline H&
$1$&\setlength{\unitlength}{4.mm}
\begin{picture}(1,1)
\put(0,0){$\bullet$}
\put(.6,.1){\begin{scriptsize}$1$\end{scriptsize}}
\end{picture}&\\\hline
\caption{Types, subtypes and posets $\ll(\al)$ for $\ga\in\op{PGL}_3(k)$}
\end{longtable}   
\end{scriptsize}
In Tables 1,2 we have denoted by $B_2$, $B_3$ the rational Jordan blocks of res\-pective dimension 2,3, with irreducible characteristic polynomial.
Also, $\mathcal{S}_G$ denotes the set of parameters describing the subtypes of type G; these are triples $(d,e,f)$ of divisors of $q-1$ satisfying
$d\ge e\ge f>1$ and $\lcm(d,e)=\lcm(d,f)=\lcm(e,f)$.

In these tables we classify the automorphisms $\ga$ of $\op{PGL}_{N+1}(k)$ (for $N=1,\,2$) in ``types", given essentially by the different Jordan normal forms. For each type we display the Hasse diagram of the poset $\igg$ and the associated ``subtypes" determined by the different values of the exponents of the nodes of this poset, indicated in the fourth column. The nodes of degree one are represented by $\bullet$ and labelled with the value of $\exp V$. The nodes with greater degree are represented by $\circ$ and labelled with $\exp V(\dg V)$. The value of $\dm V$ is given by the 
  vertical level of the node (of height $0,1$ or $0,1,2$) inside the poset; when there is some ambiguity we write the dimension of a concrete level in the left side of the poset.  
These tables furnish an intrinsic description of the set $\mathcal{S}$ of all possible subtypes, and for each $\al\in\mathcal{S}$ they exhibit the structure of the poset $\ll(\al)$, whose Moebius function is easy to compute. This accomplishes tasks (1), (2) mentioned at the end of the last section. 

Denote by $\igk$ the subposet of $\igg$ determined by the nodes of degree one. Note that $\igk=\igg$ except for a few cases where there is only one node of degree greater than one. Therefore, the equivalence relation ``having the same subtype" can be reformulated as follows (this result is not true for $N>2$)
\begin{lem}
Let $\ga,\,\ga'\in\cc$. The following conditions are equivalent:
\begin{enumerate}
\item   $\ga$ and $\ga'$ have the same subtype.
\item   there exists a poset isomorphism $\igk\iso \ll_{\ga'}^{\op{pr}}(k)$ preserving the weight $(\dim V,\exp V)$ of each node.
\item   $\ga$ and $\ga'$ have the same cycle type as permutations of $\pr2(k)$. That is, they decompose into a product of disjoint cycles of the same length.
\end{enumerate}
\end{lem}
\begin{pf}
The equivalence of items 1 and 2 is obvious.  The equivalence of items 2 and 3 is consequence of Lemma \ref{gaorbit}.  \hfill{$\Box$}
\end{pf}
For $N=2$ we computed in \cite[sec.2]{mn} the weighted sums $N_{\al}:=\sum_{\ga\in\cc_{\al}}|\g_{\ga}|^{-1} $, where 
$\cc_{\al}$ is the subfamily of $\cc$ of all $\ga$ with a concrete cycle type. For $N=1$ the computation of the weighted sums $N_{\al}$ is easily deduced from \cite[Lemma 2.2]{lmnx}. By the above lemma the classification of $\cc$ into subtypes coincides with the classification according to the cycle type, so that our universal coefficients $M_{\al}$ coincide with the numbers $N_{\al}$. This accomplishes task (3) of the end of the last section, and this allows us to use Theorem \ref{ffinal} to obtain an explicit computation of $t_1(n)$, $t_2(n)$, $\overline{t}_1(n)$, $\overline{t}_2(n)$. More precisely, restricting our attention to $t_1(n)$, $t_2(n)$, if we split the formula of Theorem \ref{ffinal} according to the different types of the elements of $\cc$  we obtain the following results:
\begin{thm}\label{t1n}
Let $\varphi$ be Euler's totient function and $P$ any point of degree two in $\pr1(\kb)$.
The numbers $t_1(n)$ split into the sum of four terms: $$t_1(n)=t_{1,A}(n)+t_{1,B}(n)+t_{1,C}(n)+t_{1,D}(n),$$ each term having the following generating function:
$$
\sum_{n\ge0}t_{1,A}(n)x^n=\frac{1+x^2}{2(q+1)}\sum_{d|q+1,\,d>1}\varphi(d)f_{\pr1\setminus O_{\sigma}(P)}(x^d),
$$
$$
\sum_{n\ge0}t_{1,B}(n)x^n=\frac{1+x}{q}f_{\af1}(x^p),
$$
$$
\sum_{n\ge0}t_{1,C}(n)x^n=\frac{(1+x)^2}{2(q-1)}\sum_{d|q-1,\,d>1}\varphi(d)f_{\af1\setminus \{\ast\}}(x^d),
$$
$$
\sum_{n\ge0}t_{1,D}(n)x^n=\frac{1}{q(q-1)(q+1)}f_{\pr1}(x).
$$
\end{thm}
Explicit expresions for the numbers $t_1(n)$ were found in \cite{lmnx}, but this computation of their generating function is new. 
\begin{thm}\label{t2n}
Let $P$ be a point of  degree two, $Q$ a point of degree three, $P_1,\,P_2,\,P_3$ points of degree one and $L,\,L'$ two non-parallel lines of $\af2$. The numbers $t_2(n)$ split into the sum of eight terms:
$$
t_2(n)=t_A(n)+t_B(n)+t_C(n)+t_D(n)+t_E(n)+t_F(n)+t_G(n)+t_H(n),
$$  
each term having the following generating function:
$$
\sum_{n\ge0}t_A(n)x^n=\frac{1+x^3}{3(q^2+q+1)}\sum_{d|q^2+q+1,\,d>1}\varphi(d)f_{\pr2\setminus O_{\sigma}(Q)}(x^d).
$$
\begin{multline*}
\sum_{n\ge0}t_B(n)x^n=\frac{(1+x)(1+x^2)}{2(q^2-1)}\left(\sum_{d|q+1,\,d>1,\,d\mbox{\scriptsize \,odd}}\varphi(d)f_{\pr2\setminus \{P_1,\,O_{\sigma}(P)\}}(x^d)+\right.\\\left.
\sum_{d|q+1,\,e|q-1,\,d,e>1}\varphi(d)\varphi(e)\delta_{d,e}f_{\pr1\setminus O_{\sigma}(P)}(x^d)f_{\af2\setminus\{\ast\}}(x^{de})\right),
\end{multline*}
where  
$$\as{1}
\delta_{d,e}=\left\{\begin{array}{ll}
1&\quad\mbox{ if $d$ odd}\\0&\quad\mbox{ if $d$ even and }e|(q-1)/2\\2&\quad\mbox{ if $d$ even and }e\nmid(q-1)/2.\\\end{array}\right.
$$
According to $p>2$ or $p=2$ we have respectively
$$
\sum_{n\ge0}t_C(n)x^n=\frac{1+x}{q^2}f_{\pr2\setminus \{\ast\}}(x^p), \ \mbox{ or }\ 
\sum_{n\ge0}t_C(n)x^n=\frac{1+x}{q^2}f_{\af1}(x^2)f_{\af2}(x^4).
$$
$$
\sum_{n\ge0}t_D(n)x^n=\frac{1}{q^3(q-1)}f_{\pr1}(x)f_{\af2}(x^p).
$$
$$
\sum_{n\ge0}t_E(n)x^n=\frac{(1+x)^2f_{\af1}(x^p)}{q(q-1)}\sum_{d|q-1,\,d>1}\varphi(d)f_{\af1\setminus \{\ast\}}(x^d)f_{\af2\setminus L}(x^{pd}).
$$
$$
\sum_{n\ge0}t_F(n)x^n=\frac{(1+x)f_{\pr1}(x)}{q(q-1)^2(q+1)}\sum_{d|q-1,\,d>1}\varphi(d)f_{\af2\setminus \{\ast\}}(x^d).
$$
\begin{multline*}
\sum_{n\ge0}t_G(n)x^n=\frac{(1+x)^3}{(q-1)^2}\left(\sum_{(d,e,f)\in\mathcal{S}_G,\,d=e=f}\frac16\varphi(m)\psi(m)f_{\pr2\setminus
\{P_1,P_2,P_3\}}(x^m)+\right.\\
+\sum_{(d,e,f)\in\mathcal{S}_G,\,d=e>f}\frac12\varphi(m)\varphi(h)\psi(H)f_{\af1\setminus \{\ast\}}(x^f)f_{\af2\setminus \{\ast\}}(x^m)+\\
\left.\sum_{(d,e,f)\in\mathcal{S}_G,\,d>e>f}\!\!\!\!\!\!\varphi(m)\varphi(h)\psi(H)f_{\af1\setminus \{\ast\}}(x^d)f_{\af1\setminus \{\ast\}}(x^e)f_{\af1\setminus \{\ast\}}(x^f)f_{\af2\setminus (L\cup L')}(x^m)\right)
\end{multline*}
where $m=\lcm(d,e)$, $(def)/m^2=hH$ is the unique decomposition of this divisor of $m$ into a product of positive divisors $h,\,H$ satisfying respectively
$$
v_{\ell}(h)<v_{\ell}(m), \quad v_{\ell}(H)=v_{\ell}(m),
$$for any prime divisor $\ell$ of $(def)/m^2$, and $\psi$ is the multiplicative function determined by $\psi(\ell^r)=(\ell-2)\ell^{r-1}$ for any prime power. Finally,
$$
\sum_{n\ge0}t_H(n)x^n=\frac{f_{\pr2}(x)}{q^3(q^2+q+1)(q-1)^2(q+1)}.
$$
\end{thm}
In these formulas we find several functions $f_{V}(x)=\op{Z}(V/k,x)/\op{Z}(V/k,x^2)$ (cf. Theorem \ref{tu})
for certain locally closed subvarieties $V$ of $\pr1$ and $\pr2$. As mentioned in the last section, one can express all these functions in terms of $f_{\pr{r}}(x)$, $r=0,1,2$, by using the Moebius function of certain poset. However, in our lower dimension cases, in order to find a concrete expression for $t_1(n)$, $t_2(n)$ it is better to use the explicit computations of the coefficients of $f_{V}(x)$ that we found in section \ref{secu} for these particular quasiprojective varieties $V$.

The computation of the generating function of $\overline{t}_1(n)$, $\overline{t}_2(n)$ is given by the same formulas, substituting the functions $f_{V}(x)$ by the corresponding functions $\overline{f}_{V}(x)=\op{Z}(V/k,x)$.

We end with a remark that might lead to further work on this topic. We showed in \cite{mn} that the numbers of $\op{PGL}_3(k)$-orbits of pointwise rational $n$-sets of the plane can be expressed as a polynomial in $q$ with rational coefficients. In contrast to this situation, Theorems \ref{t1n}, \ref{t2n} and the computations of section \ref{secu} provide formulas for $t_1(n)$, $\overline{t}_1(n)$, $t_2(n)$, $\overline{t}_2(n)$ as a polynomial in $q$ with {\it integer} coefficients. For instance, with some extra work for the computation of $t_G(7)$, our results provide the following expression for $t_2(7)$:
\begin{multline*}
t_2(7)=q^6+q^5+3q^4+6q^3+11q^2+4q+\\\\+[2(q^2+6q+9)]_{3|q-1}+[4]_{5|q+1}+[20]_{5|q-1}+[4q+6]_{4|q-1}+[8]_{7|q-1}+[6]_{7|q+1}+\\\\+
[4]_{12|q-1}+[2]_{7|q^2+q+1}+[q^2+4q-1]_{3|q}+[2]_{3|q,\,4|q-1}+[3]_{5|q}+[2]_{7|q},
\end{multline*}
if $p>2$, whereas for $p=2$ we have:
\begin{multline*}t_2(7)=
q^6+q^5+3q^4+6q^3+8q^2-1+\\\\+[2(q^2+6q+4)]_{3|q-1}+[4]_{5|q+1}+[20]_{5|q-1}+[6]_{7|q-1}.
\end{multline*}
One may speculate if the numbers $t_N(n)$, $\overline{t}_N(n)$ have this property for all $N,\,n$.

\end{document}